\documentclass[12pt]{article}

\usepackage{amsthm,amstext, amsmath,latexsym,amsbsy,amssymb}
\newtheorem{definition}{Definition}[section]

\newtheorem{lemma}{Lemma}[section] 

\numberwithin{equation}{section}
\newtheorem{theorem}{Theorem}[section] 

\newcommand{\h}{\hspace*{.24in}}

\allowdisplaybreaks
\begin{document}
\title{Parallel Schwarz Waveform Relaxation Algorithm for an N-Dimensional Semilinear Heat Equation}
\author{Minh-Binh TRAN\\
Laboratoire Analyse G\'eom\'etrie et Applications\\
Institut Galil\'ee, Universit\'e Paris 13, France\\
Email: binh@math.univ-paris13.fr}
\maketitle
      
\begin{abstract}
We present in this paper a proof of well-posedness and convergence for the parallel Schwarz Waveform Relaxation Algorithm adapted to an N-dimensional semilinear heat equation. Since the equation we study is an evolution one, each subproblem at each step has its own local existence time, we then determine a common existence time for every problem in any subdomain at any step. We also introduce a new technique: Exponential Decay Error Estimates, to prove the convergence of the Schwarz Methods, with multisubdomains, and then apply it to our problem.
\end{abstract}

\section{Introduction}
In the pioneer work \cite{Lions:1987:OSA}, \cite{Lions:1989:OSA}, \cite{Lions:1990:OSA}, P. L. Lions laid the foundations of the modern theory of Schwarz algorithms. He also proposed to use the Schwarz alternating method for evolution equations, and studied the algorithm for nonlinear monotone problems. Later, Schwarz waveform relaxation algorithms, by refering to the paper \cite{Burrage:1996:PPW}, were designed independently in \cite{Gander:1998:STC} and \cite{Giladi:2002:STD} for the linear advection-diffusion equation. They try to solve, on a given time interval, a sequence of Cauchy Problems with the transmission conditions of Cauchy type on overlapping subdomains. The algorithm is well-posed with some compatibility conditions. 
\\\h An extension to the nonlinear reaction-diffusion equation in dimension 1 was considered in \cite{Gander:1999:WRA}. For  nonlinear problems, especially evolutional equations, there are some cases that the solutions are blowed up in finite time, which means that if we devide the domain into several subdomains, at each step we can get different existence times in different domains, and we do not know if there exists a common existence time for all iterations. However, with the hypothese $f'(c)\leq C$ in \cite{Gander:1999:WRA}, we do not encounter this difficulty and the iterations are defined naturally on an unbounded time interval. Proofs of linear convergence on unbounded time domains, and superlinear convergence on finite time intervals were then given in case of n subdomains, based on some explicit computations on the linearized equations. Another extension to monotone nonlinear PDEs in higher dimensions were considered by Lui in \cite{Lui:2001:OSM}, \cite{Lui:2003:OMI}. In these papers, some monotone iterations for Schwarz methods are defined in order to get the convergence of the algorithm, based on the idea of the sub-super solutions method in PDEs and no explosion is considered. Recently, an extension to Systems of Semilinear Reaction-Diffusion Equations was investigated in \cite{Descombes:2010:SWR}. Some systems in dimension 1 were considered and the proofs of the well-posedness and the convergence of the algorithm used in this paper were a development of the technique used in \cite{Gander:1999:WRA}.
\\\h We consider here the semilinear heat equation $(\ref{2e1a})$, in a domain $\Omega=D\times(a,b)$ of $\mathbb{R}^N$, with the nonlinearity of the form $(\ref{2e1b})$, which allows explosion of solutions in finite time. We cut the domain into bands $\Omega_i=D\times(a_i,b_i)$, with $a_1=a$ and $b_I=b$. These bands are ovelapping, i.e. for all $i\in\{1,I-1\}$, $a_{i+1}<b_i<b_{i+1}$. In each of these subdomains, we solve a heat equation with Dirichlet limit conditions. Since the domains are not smooth, we cannot use classical results about semilinear heat equations on smooth domains; we then establish some new proofs of existence for a general domain in Theorem 2.1. Applying the results in Theorem 2.1 for the equation $(\ref{2e2})$, we get an existence result for a semilinear heat equation in a domain of the type $\Omega=D\times(a,b)$ in Theorem 2.2.
\\\h Theorem 2.3 confirms that the algorithm is well-posed and there exists a common existence time for all subdomains at all iterations despite of the phenomenon of explosion. The common existence time $T^*$ is computed explicitely so that one can use it in numerical simulations. This is a collolary of Theorem 2.2.
\\\h We prove in Theorem 2.4 that the algorithm converges linearly. There are five main techniques to prove the convergence of Domain Decomposition Algorithms and they are: Orthorgonal Projection used for a linear Laplace equation (see \cite{Lions:1987:OSA}), Fourier and Laplace Transforms used for linear equations (see  \cite{Gander:2007:OSW},\cite{Gander:1998:OCO}, \cite{Gander:1999:OSM}, \cite{Giladi:2002:STD}), Maximum Priciple used for linear equations (see \cite{Gander:2002:OSW}), Energy Estimates used for nonoverlapping algorithms (see \cite{Benamou:1997:DDM}), and Monotone Iterations (see \cite{Lui:2001:OSM}, \cite{Lui:2003:OMI}) used for nonlinear monotone problems. The convergence problem of overlapping algorithms for nonlinear equations is still open up to now. In this paper, we introduce a new technique: Exponential Decay Error Estimates, basing on the idea of constructing some Controlling Functions, that allows us to prove that the algorithm converges linearly.
\section{The main results}
We consider the semilinear heat equation
\begin{equation}\label{2e1a}
\partial_t u -\Delta u - f(u) =0,
\end{equation}
with the assumptions on $f$:
\begin{equation}\label{2e1b}
\mbox{ $f$ is in  $C^1(\mathbb{R})$ and there exists $C_f>0$, $p>1$ such that $|f'(x)|\leq C_f|x|^{p-1}$}.
\end{equation}
We first set an existence theorem for the initial boundary value problem, and more important, new estimates on the solution.
We need here some notations. We set $p_1=\frac{3(p-1)}{4p}$, $\alpha=\frac{1}{2}(\frac{1}{p-1}-\frac{3}{4p})$.
$l_1$ and  $l_2$ are positive numbers such that $\frac{1}{l_1}+\frac{1}{l_2}=1$ and $l_1 p_1 < 1$. We denote by $||u||_{k,h}$ the norm $||u||_{L^{k}(0,T,L^h(\omega))}$, when $u$ belongs to $L^{k}(0,T,L^h(\omega))$, where $\omega$ is some domain in $\mathbb{R}^N$, $k,h$ can be $\infty$. We define
\begin{equation}\label{2e1c}
\begin{array}{c}
\mbox{\large{$\tau$}}(r,m)=
 [(4\pi)^{-p_1}\ \frac{2^{p_1+p\alpha}}{1-(p_1+p\alpha)}\  C_f\max(1,2^{p-2}) (4r+ 4m)\  ]
^{-\frac{8p}{3+p}},\\[2mm]
G(r; T,m_1,m_2)=
\left(
\frac{(4\pi)^{-p_1l_1}T^{1-p_1l_1}}
    {1-p_1l_1}
\right)^{-\frac{l_2}{l_1}}
\int_0^r
\left[C_f\max\{1,2^{p-2}\}(m_1
+\zeta^{\frac{p-1}{l_2}})\zeta^{\frac{1}{l_2}}
+m_2\right]^{-l_2}
 d\zeta.
\end{array}
\end{equation}
Consider the problem
\begin{equation}\label{2e1}
\left\{
\begin{array}{ll}\partial_tw-\Delta w=f(w+v)
&\mbox{ in }{\mathcal O}\times(0,T),\\
w =0&\mbox{ on }\partial{\mathcal O}\times[0,T],\\
w(.,0)=0&\mbox{ in }\mathcal O.
\end{array}\right.
\end{equation}
\begin{theorem}\label{2t1}
Let ${\mathcal O}$ be a bounded domain in $\mathbb{R}^N$, $m({\mathcal O})$ is its measure. 
\\ Suppose $v \in  C([0,T_0,L^{2}](\Omega))$ and $|v|\leq M$ a.e.. Denote $R_1=$ $2\max_{|\zeta|\leq M}$ $|f(\zeta)|m({\mathcal O})^{\frac{1}{2}}$ $\frac{T_0^{\frac{p+3}{4p}}}{\frac{3(p-1)}{4p}}$. Then, there exists a local time
$T_*= \min (T_0,\,1,\,
\mbox{\large{$\tau$}}(R_1,M^{p-1}m({\cal O})^{\frac{p-1}{2p}}))$,
such that for all $T<T_*$, equation ($\ref{2e1}$) has a unique solution $w$ in $L^{\infty}({\cal O}\times(0,T))\cap C([0,T],L^2({\cal O}))\cap L^{2}(0,T,H_0^1({\cal O}))$ and $\partial_t w\in L^{2}(0,T,L^2({\cal O}))$.
Furthermore,
$||w||_{\infty,\infty}\leq M_*$, where
\begin{equation}\label{2e1d}
\begin{array}{c}
M_* := \bigl( \frac{4T_*}{\pi^3}\bigr)^{\frac{1}{4}}\
\ [\ C_f\ \max(1,2^{p-2})\
(G^{-1}(T_*)^{\frac{p-1}{l_2}}
\ +M^{p-1}m({\mathcal O})^{\frac{p-1}{2p}})\ G^{-1}(T_*)^{\frac{p-1}{l_2}}\\
\hspace{20mm}
+m({\mathcal O})^{\frac{1}{2}}\max_{|\zeta|\le M}|f(\zeta)|\,]\\
G^{-1}(T_*)\equiv G^{-1}(T_*; T,M^{p-1}m({\mathcal O})^{\frac{p-1}{2p}},m({\mathcal O})^{\frac{1}{2}}\max_{|\zeta|\le M}|f(\zeta)| ).
\end{array}
\end{equation}
\\ Suppose $v \in  C([0,T_0],L^{2}(\Omega))\cap L^{\infty}(0,T_0,L^{2p}(\Omega))$ a.e.. Denote $R_2=2\frac{T_0^{\frac{p+3}{4p}}}{\frac{3(p-1)}{4p}}||f(v)||_{\infty,2}$. Then, there exists a local time
$T_*= \min (T_0,\,1,\,
\mbox{\large{$\tau$}}(R_2,||v||^{p-1}_{L^{\infty}(0,T,L^{2p}(\Omega)})$,
such that for all $T<T_*$, equation $(\ref{2e1})$ has a unique solution $w$ in $L^{\infty}({\mathcal O}\times(0,T))\cap C([0,T],L^2({\mathcal O}))\cap L^{2}(0,T,H_0^1({\mathcal O}))$ and $\partial_t w\in L^{2}(0,T,L^2({\mathcal O}))$.
Furthermore,
$||w||_{\infty,\infty}\leq M_*$, where
\begin{equation}\label{2e1e}
\begin{array}{c}
M_* := \bigl( \frac{4T_*}{\pi^3}\bigr)^{\frac{1}{4}}\
\ [\ C_f\ \max(1,2^{p-2})\
(G^{-1}(T_*)^{\frac{p-1}{l_2}}
\ +||v||^{p-1}_{\infty,2p})\ G^{-1}(T_*)^{\frac{p-1}{l_2}}\\
\hspace{20mm}
+||f(v)||_{\infty,2}\,]\\
G^{-1}(T_*)\equiv G^{-1}(T_*; T,||v||^{p-1}_{\infty,2p},||f(v)||_{\infty,2})
\end{array}
\end{equation}
\end{theorem}
We consider now a bounded domain of the following form $\Omega=D\times(a,b) \subset \mathbb{R}^N$, where $D$ is a bounded domain with smooth enough boundary $\partial D$ in $\mathbb{R}^{N-1}$.  The boundary $\partial \Omega$ of $\Omega$ is made of three parts,
$\Gamma_L=\bar D\times \{a\}$, $\Gamma_R=\bar D\times \{b\}$, and $\Gamma_C=\partial D \times (a,b)$.
Dirichlet data $g$ are given on the boundary $\partial \Omega\times (0,T)$, defined by
$g_L $ on $\Gamma_L$, $g_R$ on $\Gamma_R$, $g_C$ on $\Gamma_C$. These functions are all continuous.
We now introduce the basic initial boundary value problem for $(\ref{2e1})$:
\begin{equation}\label{2e2}
\begin{array}{ll}
\partial_t u -\Delta u =f(u) &\mbox{ in }\Omega\times(0,T),\\
u =g  &\mbox{ in }\partial\Omega \times(0,T),\\
u(.,0)=u_0 &\mbox{ in }\Omega.
\end{array}
\end{equation}

\begin{theorem}\label{2t2}
Let $u_0$ in $C(\overline{\Omega})$ and $g$ in $C(\overline{\partial \Omega})$ with ${u_0}_{\,\mid\partial \Omega} =g_{\,\mid t=0}$. Let $M$ be a positive constant such that
$M>\max (||u_0||_{\infty}, ||g||_{\infty})$. Let $R_1$ be like in Theorem $\ref{2t1}$. Then, there exists a limit time $T_*= \min (1,\,$ $
\mbox{\large{$\tau$}}(R,M^{p-1}m(\Omega)^{\frac{p-1}{2p}}))$,
such that for all $T<T_*$, equation $(\ref{2e2})$ has a solution $u$ in  $L^{\infty}(\Omega\times(0,T))\cap C([0,T],L^2(\Omega))\cap L^{2}(0,T,H_0^1(\Omega))$ and $\partial_t u\in L^{2}(0,T,L^2(\Omega))$. Furthermore $u$ is continuous on $\overline{\Omega\times(0,T)}$ and $||u||_{\infty,\infty}$ $\leq M+M^*$
where $M^*$ is obtained from $(\ref{2e1d})$ by replacing ${\mathcal O}$ by $\Omega$.
\end{theorem}

We divide the domain $\Omega$ into $I$ subdomains with $\Omega_i=D\times(a_i,b_i)$, with $a_1=a$ and $b_I=b$.
We suppose for each $i\in\{1,\dots,I-1\}$, $a_{i+1}<b_i<b_{i+1}$, we denote by $L_i$ the length of $\Omega_i$: $L_i=b_{i}-a_{i}$, and by $S_i$ the size of the overlap $S_i=b_i-a_{i+1}$. 
\\ The parallel  Schwarz Waveform Relaxation Algorithm solves $I$ equations in $I$ subdomains instead of solving directly the main problem $(\ref{2e2})$. The iterate $\#k$ in the $j$-th domain, denoted by $u_{j}^{k}$, is defined by
\begin{equation}\label{2e4}
\left \{
\begin{array}{ll}
\partial_tu_{j}^{k}-\Delta u_{j}^{k}=f(u_{j}^{k})&\mbox{ in }\Omega_j\times(0,T),\vspace{.1in}\\
u_{j}^{k}(\cdot,a_j,\cdot)=u_{j-1}^{k-1}(\cdot,a_j,\cdot) &\mbox{ in }D\times(0,T),\vspace{.1in}\\
u_{j}^{k}(\cdot,b_j,\cdot)=u_{j+1}^{k-1}(\cdot,b_j,\cdot) &\mbox{ in }D\times(0,T).\vspace{.1in}
\end{array}
\right.
\end{equation}
Each iterate inherits the boundary conditions and the initial values of $u$:
\[
u_j^k= g \mbox{ on }\partial\Omega_j\cap\partial\Omega \times(0,T),\quad
u_j^k(\cdot,\cdot,0)= u_0 \mbox{ in } \Omega_j,
\]
which imposes a special treatment for the extreme subdomains,
\[
u_1^k(\cdot,a,\cdot)=g,
\quad
u_I^k(\cdot,b,\cdot)=g.
\]
An initial guess is provided, \textit{i.e.} we solve at step $0$ equations $(\ref{2e4})$, with boundary data on left and right
\[
u_{j}^{0}(\cdot,a_j,\cdot)=g_{j}^{0} \mbox{ in }D\times(0,T),\vspace{.1in}\\
u_{j}^{k}(\cdot,b_j,\cdot)=h_{j}^{0} \mbox{ on }D\times(0,T).\vspace{.1in}
\]
\begin{definition}\label{2d1}
The parallel  Schwarz waveform relaxation algorithm is well-posed  if there exists a local time $T^* \le T_*$ such that for all $T<T^*$, each subproblem $(\ref{2e4})$ in each iteration has a solution over the time interval $(0,T)$, and the set of solutions $\{u_{j}^{k},\ j\in \overline{1,I},\ k\in N\}$ is bounded in $C(\overline\Omega\times[0,T])$.
\end{definition}

Let $M$ a positive number. According to theorem \ref{2t2}, the following problem has a solution  $\phi_M$  in some interval  $C(\overline{\Omega}\times[0,T_0])$:
\begin{equation}\label{2e5}
    \left \{
        \begin{array}{ll}
            \partial_t\phi_M-\Delta\phi_M=f(\phi_M)
                &\mbox{ in }\Omega\times(0,T_0),\\
            \phi_M =M
                &\mbox{ in }\partial\Omega\times[0,T_0],\\
            \phi_M(\cdot,\cdot,0)=M
                &\mbox{ in }\bar\Omega.
        \end{array}
    \right.
\end{equation}

The next theorem gives a common existence time for the iterates:

\begin{theorem} \label{2t3}
Let $M$ be a positive constant such that
\[
M>\max (||u_0||_{\infty},\  ||g||_{\infty},\ (||g_i^0||_{\infty},\ ||h_i^0||_{\infty}, i\in \overline{1,I})\ ).
\]
Suppose the data $u_0$ and $g_j$, $g_j^0$ and $h_j^0$ are continuous and satisfy the compatibility conditions
\[
{u_0}_{\,\mid\partial \Omega} =g_{\,\mid t=0},\
u_0(\cdot,a_j)={g_j^0}_{\,|\mid t=0},\
u_0(\cdot,b_j)={h_j^0}_{\,\mid t=0}.\
\]
Let $M_*$ be greater than the maximum of $\phi_M$ on the boundaries of $\Omega\times(0,T)$ and $\Omega_j\times(0,T)$, $j\in \overline{1,I}$.
Put
$T^*= \min (T_0,\,1,\,T_*,\mbox{\large{$\tau$}}(R_1,M_*^{p-1}m({\Omega_j})^{\frac{p-1}{2p}}), j\in\overline{1,I}
))$.
Then the parallel  Schwarz waveform relaxation algorithm $(\ref{2e4})$ is well-posed with a local time at least equal to $T^*$.
\end{theorem}

We denote by $e_j^k$ the error $u_j^k-u$, where $u$ is the solution of $(\ref{2e2})$. Let $P$ be a function from $\mathbb{R}$ to $\mathbb{R}$ such that (i) $P\in  C^2(\mathbb{R})$; (i) $P(x),P''(x)\geq 0$ $\forall x\in\mathbb{R}$; (iii) $P(x)=0$ iff $x=0$ and $P'(0)=0$; (iv) $\forall M>0$, there exists $K(M)$ such that $\left|\frac{xP'(x)}{P(x)}\right|<K(M)$, $\forall$ $x\in\mathbb{R}$, $|x|<M$. We finally state the convergence of the algorithm:

\begin{theorem} \label{2t4}
With the same assumptions as in Theorem \ref{2t2}, let $\gamma$ be a constant large enough and denote by $\bar\epsilon$ the constant $\sqrt \frac{\gamma}{2}\frac{S_1\dots S_{I-1}}{L_2\dots L_{I-1}}$. If we put $E_k=$ $\max_{j\in\overline{1,I}}$ $||P(e_j^k)\exp(-\gamma t)||_{C(\overline{\Omega_j\times(0,T)})}$ $=$ $\max_{j\in\overline{1,I}}$ $||P(u_{j}^{k}-u)\exp(-\gamma t)||_{C(\overline{\Omega_j\times(0,T)})}$ and $\bar E_k=\max_{j\in\overline{1,I-1}}\{E_{k+j}\},$
the parallel  Schwarz waveform relaxation algorithm $(\ref{2e4})$ converges linearly in the following sense: for every $T < T^*$,
$$\bar E_{n}\leq \bar E_0\exp(-n\bar\epsilon),\forall n\in\mathbb{N},$$
and as a consequence
\[
\lim_{k\to\infty}
\max_{j\in[1,I]}
\|u_{j}^{k}-u\|_{C(\overline{\Omega_j}\times[0,T])}=0.
\]
\end{theorem}
\section{Proof of The Existence Results for Semilinear Heat Equation in an Arbitrary Bounded Domain - Theorem $\ref{2t1}$} 
\subsection{Premiliary Results} Let $\omega$ be an arbitrary bounded domain in $\mathbb{R}^N$ and consider the following equation
\begin{equation}
\left \{ \begin{array}{ll}\partial_t\zeta-\Delta\zeta=0&\mbox{ in }\omega\times(0,T),\vspace{.1in}\\ 
\zeta(.,.)=0&\mbox{ on }\partial\omega\times[0,T],\vspace{.1in}\\ 
\zeta(.,0)=\zeta_0&\mbox{ on }\bar\omega.\end{array}\right. 
\label{3e1}
\end{equation}
\h We consider the operator $A\zeta=-\Delta\zeta$, $D(A)=\{\zeta|\zeta\in H_0^1(\omega),\Delta\zeta\in L^2(\omega)\}$. According to Proposition 2.6.1 page 26 $\cite{Cazenave:1998:AIS}$, $A$ is an m-dissipative operator with dense domain in $L^2(\omega)$. Let $S(t)$ be the Dirichlet semigroup associated with $A$ on $L^2(\omega)$ (see \cite{Davis:1989:HKS}). If $\zeta_0\in D(A)$, due to Theorem 3.1.1 page 33 \cite{Cazenave:1998:AIS}, $\zeta(t)=S(t)\zeta_0$ is a solution of $(\ref{3e1})$. 
\\\h According to Section 6.5 page 334 \cite{Evans:1998:PDE}, there exist a sequence of eigenvalues $0<\lambda_1\leq\lambda_2\dots$ and $\lambda_k\to\infty$ as $k\to\infty$ and an othornormal basis $\{w_k\}_{k=1}^{\infty}$ of $L^2(\omega)$ which is also an orthogornal basis of $H_0^1(\omega)$, where $w_k$ is the eigenfunction corressponding to $\lambda_k$:
\begin{eqnarray*}
\left \{ \begin{array}{ll}Aw_k=\lambda_kw_k&\mbox{ in }\omega,\vspace{.1in}\\ 
w_k=0&\mbox{ on }\partial\omega.\end{array}\right. 
\end{eqnarray*}
\h Denote by $<,>$ the scalar product in $L^2(\omega)$, we have the following Lemma
\begin{lemma}
If $\zeta_0\in L^2(\omega)$, then we have that 
\begin{equation}
S(t)\zeta_0=\sum_{k=1}^{\infty}e^{-t\lambda_k}<\zeta_0,w_k>w_k.
\label{3e2}
\end{equation}
\label{3l1}
\end{lemma}
\begin{proof} We prove the Lemma in two steps.
{\\\h\bf Step 1:} Suppose that $\zeta_0\in D(A)$ and put
$$\nu(x,t)=\sum_{k=1}^{\infty}e^{-t\lambda_k}<\zeta_0,w_k>w_k,$$
we prove that $\nu$ is a solution of $(\ref{3e1})$.
\\\h We compute the norm of  $\nu(.,t)$ in $L^{2}(\omega)$
$$||\nu(x,t)||_{L^{2}(\omega)}=\sum_{k=1}^{\infty}(e^{-t\lambda_k}<\zeta_0,\omega_k>)^2\leq \sum_{k=1}^{\infty}(<\zeta_0,\omega_k>)^2=||\zeta_0||^2_{L^2(\omega)}.$$
 The inequality implies $\nu\in L^{\infty}(0,T,L^2(\omega))$. 
\\\h We now need that $\nu\in L^2(0,T,H_0^1(\omega))$, and it is enough to prove
$$\sum_{k=1}^{\infty}\int_0^T(e^{-t\lambda_k}<\zeta_0,w_k>)^2||w_k||^2_{H_0^1(\omega)}dt<\infty.$$
 We estimate the term on the left hand side of the previous inequality
\begin{eqnarray*}
\sum_{k=1}^{\infty}\int_0^T(e^{-t\lambda_k}<\zeta_0,w_k>)^2||w_k||^2_{H_0^1(\omega)}dt & = &\sum_{k=1}^{\infty}\int_0^Te^{-2t\lambda_k}\lambda_kdt<\zeta_0,w_k>^2\\
& = & \sum_{k=1}^{\infty}\frac{1}{2}(1-e^{-2T\lambda_k})<\zeta_0,w_k>^2\\
& \leq & \frac{1}{2}||\zeta_0||^2_{L^2(\omega)}.
\end{eqnarray*}
 The estimate leads to $\nu\in L^2(0,T,H_0^1(\omega))$.
\\\h Since $\nu(0,.)=\sum_{k=1}^{\infty}<\zeta_0,w_k>w_k=\zeta_0,$ we only have to prove that $\partial_t\nu-\Delta\nu=0$. In order to do this, we put $\nu_n=\sum_{k=1}^n(e^{-t\lambda_k}<\zeta_0,w_k>)w_k$. We compute by a direct manner 
$$\partial_t\nu_n=\sum_{k=1}^n-\lambda_k(e^{-t\lambda_k}<\zeta_0,w_k>)w_k=\sum_{k=1}^n(e^{-t\lambda_k}<\zeta_0,w_k>)\Delta w_k=\Delta\nu_n.$$
It follows from the previous equation that
\begin{eqnarray*}
\int_0^T\int_{\omega}(\partial_t\nu_n)^2dxdt & = &-\int_0^T\int_{\omega}\nabla\nu_n\nabla(\partial_t\nu_n)dxdt\\
& = &-\frac{1}{2}\int_{\omega}(\nabla\nu_n)^2(T)dx+\frac{1}{2}\int_{\omega}(\nabla\nu_n)^2(0)dx.
\end{eqnarray*}
 The inequality $\int_{\omega}(\nabla\nu_n)^2(0)dx=\int_{\omega}\sum_{k=1}^n<\zeta_0,w_k>^2(\nabla w_k)^2dx\leq$$ \int_{\omega}\sum_{k=1}^{\infty}<\zeta_0,w_k>^2(\nabla w_k)^2dx=\int_{\omega}(\nabla\zeta_0)^2dx=||\zeta_0||^2_{H_0^1(\omega)}$ implies $\int_0^T\int_{\omega}(\partial_t\nu_n)^2dxdt\leq||\zeta_0||^2_{H_0^1(\omega)}$. We obtain $\partial_t\nu\in L^2(0,T,L^2(\omega))$ and then $\partial_t\nu=\Delta\nu$ in $L^2(0,T,L^2(\omega))$.
{\\\h\bf Step 2:} We prove $(\ref{3e2})$.
\\\h We now put $\mu=\nu-\zeta$ and recall that $\zeta=S(t)\zeta_0$ is a solution of $(\ref{3e1})$. It follows that $\mu$ is the solution of the following equation
\begin{equation}
\left \{ \begin{array}{ll}\partial_t\mu-\Delta\mu=0&\mbox{ in }\omega\times(0,T),\vspace{.1in}\\ 
\mu(.,.)=0&\mbox{ on }\partial\omega\times[0,T],\vspace{.1in}\\ 
\mu(.,0)=0&\mbox{ on }\bar\omega,\vspace{.1in}\\
\mu\in L^2(0,T,H_0^1(\omega))$, $\partial_t\mu\in L^2(0,T,L^2(\omega)).
\end{array}\right. 
\label{3e3}
\end{equation}
We derive from $(\ref{3e3})$ that
$$\int_0^t<\partial_t\mu(.,s),\mu(.,s)>ds+\int_0^t<\nabla\mu(.,s),\nabla\mu(.,s)>ds=0.$$
This means $\frac{1}{2}||\mu(t)||^2_{L^2(\omega)}+\int_0^t||\nabla\mu||^2_{L^2(\omega)}=0.$ Then $\mu=0$, or $\nu=\zeta$; and $(\ref{3e2})$ holds for $\zeta_0\in D(A)$. Since $D(A)$ is dense in $L^2(\omega)$ and $S(t)\in{\cal{L}}(L^2(\omega))$, $(\ref{3e2})$ holds for $\zeta_0\in L^2(\omega)$.
\end{proof}
\begin{lemma} Suppose that $g\in C([0,T],L^2(\omega))$ and $\rho_0\in D(A)$, we consider the following equation
\begin{equation}
\left \{ \begin{array}{ll}\partial_t\rho-\Delta\rho=g&\mbox{ in }\omega\times(0,T),\vspace{.1in}\\ 
\rho(.,.)=0&\mbox{ on }\partial\omega\times[0,T],\vspace{.1in}\\ 
\rho(.,0)=\rho_0&\mbox{ on }\bar\omega.\end{array}\right. 
\label{3e4}
\end{equation}
\h Then $(\ref{3e4})$ has a solution in $L^{2}(0,T,H_0^1(\omega))\cap L^{\infty}(0,T,L^2(\omega))$ given by the following formula
\begin{equation}
\rho(.,t)=S(t)\rho_0+\int_0^tS(t-s)g(.,s)ds.
\label{3e5}
\end{equation}
\label{3l2}
\end{lemma}
\begin{proof} According to Lemma $\ref{3l1}$, $S(t)\rho_0$ is the solution of $(\ref{3e1})$ with the initial condition $\rho_0$ and $S(t)\rho_0\in L^{2}(0,T,H_0^1(\omega))\cap L^{\infty}(0,T,L^2(\omega))$. We can suppose with out loss of generality that $\rho_0=0$ in our proof. Then, by Lemma $\ref{3l1}$, formula $(\ref{3e5})$ can be rewritten in the following form
\begin{equation}
\rho(.,t)=\sum_{k=1}^{\infty}\left[\int_0^te^{-(t-s)\lambda_k}<g(.,s),w_k>ds\right]w_k.
\label{3e6}
\end{equation}
{\\\h\bf Step 1:} We prove that $\rho\in L^2(0,T,H_0^1(\omega))\cap L^{\infty}(0,T,L^2(\omega))$.
\\\h The fact that $\rho\in L^{\infty}(0,T,L^2(\omega))$ can be derived from the following inequality
\begin{eqnarray*}
\sum_{k=1}^{\infty}\left[\int_0^te^{-(t-s)\lambda_k}<g(.,s),w_k>ds\right]^2 & \leq &\sum_{k=1}^{\infty}\int_0^te^{-2(t-s)\lambda_k}ds\int_0^t<g,w_k>^2ds\\
& \leq & T||g||^2_{L^2(0,T,L^2(\omega))}.
\end{eqnarray*}
\\\h Moreover, from
\begin{eqnarray*}
||\rho(t)||^2_{H_0^1(\omega)}& = &\sum_{k=1}^{\infty}\left[\int_0^te^{-(t-s)\lambda_k}<g(.,s),w_k>ds\right]^2||w_k||^2_{H_0^1(\omega)}\\
& = &\sum_{k=1}^{\infty}\left[\int_0^te^{-(t-s)\lambda_k}<g(.,s),w_k>ds\right]^2\lambda_k\\
& \leq &\sum_{k=1}^{\infty}\lambda_k\int_0^te^{-2(t-s)\lambda_k}ds\int_0^t<g,w_k>^2ds\\
& \leq & \frac{1}{2}||g||^2_{L^2(0,T,L^2(\omega))},
\end{eqnarray*}
we obtain that $\rho\in L^{2}(0,T,H_0^1(\omega))$.
{\\\h\bf Step 2:} We prove that $\partial_t\rho-\Delta\rho=g$.
\\\h Put $\rho_n(t)=\sum_{k=1}^n[\int_0^te^{-(t-s)\lambda_k}<g(.,s),w_k>ds]w_k$ and $g_n(t)=\sum_{k=1}^n<g(t),w_k>w_k$. Since $g\in C([0,T],L^2(\omega))$, $\partial_t\rho_n(t)=\Delta\rho_n(t)+g_n(t).$ This formula implies
$$\int_0^T\int_{\omega}(\partial_t\rho_n)^2dxdt=-\int_{\omega}(\nabla\rho_n)^2(T)dx+\int_{\omega}(\nabla\rho_n)^2(0)dx+\int_0^T\int_{\omega}g_n\partial_t\rho_ndxdt.$$
\\ We estimate the RHS of the previous equation
$$\int_0^T\int_{\omega}(\partial_t\rho_n)^2dxdt\leq\int_0^T\int_{\omega}g_n\partial_t\rho_ndxdt\leq\frac{1}{2}\int_0^T\int_{\omega}(\partial_t\rho_n)^2dxdt+\frac{1}{2}\int_0^T\int_{\omega}g_n^2dxdt,$$
\\ and deduce $\int_0^T\int_{\omega}(\partial_t\rho_n)^2dxdt\leq\int_0^T\int_{\omega}g_n^2dxdt.$ Then $\partial_t\rho\in L^2(0,T,L^2(\omega))$ and $\partial_t\rho-\Delta\rho=g$ in $L^2(0,T,L^2(\omega))$.
\end{proof}
\subsection{\bf Proof of Theorem \ref{2t1}} We consider again the operator $A\zeta=-\Delta\zeta$, $D(A)=\{\zeta|\zeta\in H_0^1({\mathcal O}),\Delta\zeta\in L^2({\mathcal O})\}$ with its associated Dirichlet semigroup $S(t)$. Thanks to Lemma 3.2, instead of considering directly the equation $(\ref{2e1})$, in some of the following steps, we will consider the following equation
\begin{equation}
u(t)=S(t)(0)+\int_0^tf(u+v)(s)ds=\int_0^tf(u+v)(s)ds.
\label{3e7}
\end{equation}
{\\\h\bf Step 1:} We use a fixed point argument to prove that $(\ref{3e7})$ has a solution.
\\\h We work with the case $v\in L^{\infty}(0,T,L^{2p}({\mathcal O}))$. The case $v\in L^{\infty}(0,T,L^{\infty}({\mathcal O}))$ can be treated with exactly the same manner.
\\\h We consider the Banach space $Y_T=\{u\in L^{\infty}_{loc}((0,T),L^{2p}({\mathcal O})),||u||_{Y_T}<\infty\}$ and $||u||_{Y_T}=\sup_{0<t<T}t^{\alpha}||u(t)||_{2p}$. Let $B$ be the closed ball in $Y_T$ with center $0$ and radius $R_2$, we will use the Banach Fixed Point Theorem for the mapping
\begin{eqnarray}\label{3e7b}
\Phi& : &B\to B\\\nonumber
\Phi(u)(t)& = &\int_0^tS(t-s)f(u+v)(s)ds.
\end{eqnarray}
\h Let $u_1$, $u_2$ be two functions in $B$, we now prove that $\Phi$ is a contraction
\begin{eqnarray*}
t^{\alpha}||\Phi(u_1)(t)-\Phi(u_2)(t)||_{2p}\leq t^\alpha\int_0^t||S(t-s)(f(u_1+v)-f(u_2+v))||_{2p}ds.
\end{eqnarray*}
Using the $L^p-L^q$ estimate (see Proposition 48.4 page 441 \cite{Quittner:2007:SPP}), we obtain
\begin{eqnarray*}
||S(t-s)(f(u_1+v)(s)-f(u_2+v)(s))||_{2p}\leq(4\pi(t-s))^{-\frac{3(p-1)}{4p}}||f(u_1+v)-f(u_2+v)||_2.\end{eqnarray*}
This leads to
\begin{eqnarray*}
& &t^{\alpha}\int_0^t||S(t-s)(f(u_1+v)-f(u_2+v))||_{2p}ds\\
& \leq &t^{\alpha}(4\pi)^{-\frac{3(p-1)}{4p}}\int_0^t(t-s)^{-\frac{3(p-1)}{4p}}||f'(v+\zeta)||_{\frac{2p}{p-1}}||u_1-u_2||_{2p}ds\\
& \leq &t^{\alpha}(4\pi)^{-\frac{3(p-1)}{4p}}\int_0^t(t-s)^{-\frac{3(p-1)}{4p}}||C_f|u_1+v|^{p-1}+C_f|u_2+v|^{p-1}||_{\frac{2p}{p-1}}||u_1-u_2||_{2p}ds\\
& \leq &t^{\alpha}(4\pi)^{-\frac{3(p-1)}{4p}}C_f\int_0^t(t-s)^{-\frac{3(p-1)}{4p}}(||u_1+v||_{2p}^{p-1}+||u_2+v||_{2p}^{p-1})||u_1-u_2||_{2p}ds.
\end{eqnarray*}
We have $||u_i+v||_{2p}^{p-1}\leq(||u_i||_{2p}+||v||_{2p})^{p-1}$, for $i=1,2$. Moreover, if $p\geq 2$, we have that $(||u_i||_{2p}+||v||_{2p})^{p-1}\leq 2^{p-2}(||u_i||_{2p}^{p-1}+||v||_{2p}^{p-1})$ and if $p<2$, we have that $(||u_i||_{2p}+||v||_{2p})^{p-1}\leq ||u_i||_{2p}^{p-1}+||v||_{2p}^{p-1}$, for $i=1,2$. These inequalities lead to
\begin{eqnarray*}
& &t^{\alpha}\int_0^t||S(t-s)(f(u_1+v)(s)-f(u_2+v)(s))||_{2p}ds\\
& \leq & t^{\alpha}(4\pi)^{-\frac{3(p-1)}{4p}}C_f\max\{2^{p-2},1\}\times\\
& &\times\int_0^t(t-s)^{-\frac{3(p-1)}{4p}}(||u_1||_{2p}^{p-1}+||u_2||_{2p}^{p-1}+2||v||_{2p}^{p-1})||u_1-u_2||_{2p}ds\\
& \leq & t^{\alpha}(4\pi)^{-\frac{3(p-1)}{4p}}C_f\max\{2^{p-2},1\}\int_0^t(t-s)^{-\frac{3(p-1)}{4p}}s^{-(p-1)\alpha}\times\\
& &\times(s^{(p-1)\alpha}||u_1||_{2p}^{p-1}+s^{(p-1)\alpha}||u_2||_{2p}^{p-1}+2s^{(p-1)\alpha}||v||_{2p}^{p-1})||u_1-u_2||_{2p}ds\\
& \leq & t^{\alpha}(4\pi)^{-\frac{3(p-1)}{4p}}C_f\max\{2^{p-2},1\}\times\\
& & \times\int_0^t(t-s)^{-\frac{3(p-1)}{4p}}s^{-(p-1)\alpha}(2R_2+2T^{(p-1)\alpha}||v||_{2p}^{p-1})||u_1-u_2||_{2p}ds\\
& \leq &t^{\alpha}(4\pi)^{-\frac{3(p-1)}{4p}}C_f\max\{2^{p-2},1\}(2R_2+2T^{(p-1)\alpha}||v||_{{\infty,2p}}^{p-1})\times\\
& &\times\int_0^t(t-s)^{-\frac{3(p-1)}{4p}}s^{-(p-1)\alpha}||u_1-u_2||_{2p}ds\\
& \leq & t^{\alpha}(4\pi)^{-\frac{3(p-1)}{4p}}C_f\max\{2^{p-2},1\}(2R_2+2T^{(p-1)\alpha}||v||_{{\infty,2p}}^{p-1})\times\\
& &\times\int_0^t(t-s)^{-\frac{3(p-1)}{4p}}s^{\alpha}||u_1-u_2||_{2p}s^{-p\alpha}ds\\
& \leq & (4\pi)^{-\frac{3(p-1)}{4p}}C_f\max\{2^{p-2},1\}(2R_2+2T^{(p-1)\alpha}||v||_{{\infty,2p}}^{p-1})||u_1-u_2||_{Y_T}\times\\
& &\times t^{\alpha}\int_0^t(t-s)^{-\frac{3(p-1)}{4p}}s^{-p\alpha}ds,
\end{eqnarray*}
noticing that here we use the fact: $u_1$, $u_2$ $\in$ $B$, i.e. $s^{(p-1)\alpha}||u_1||_{2p}^{p-1}$, $s^{(p-1)\alpha}||u_2||_{2p}^{p-1}$ $\leq$ $R_2$.\\
We consider the last term on the RHS of the previous inequality
\begin{eqnarray*}
& & t^{\alpha}\int_0^t(t-s)^{-\frac{3(p-1)}{4p}}s^{-p\alpha}ds\\
& = &\int_0^tt^{1+\alpha-p\alpha-\frac{3(p-1)}{4p}}\left(1-\frac{s}{t}\right)^{-\frac{3(p-1)}{4p}}\left(\frac{s}{t}\right)^{-p\alpha}d\left(\frac{s}{t}\right)\\
& = &t^{1+\alpha-p\alpha-\frac{3(p-1)}{4p}}\int_0^1(1-\nu)^{-\frac{3(p-1)}{4p}}\nu^{-p\alpha}d\nu\\
& = &t^{1+\alpha-p\alpha-\frac{3(p-1)}{4p}}\left[\int_0^{\frac{1}{2}}(1-\nu)^{-\frac{3(p-1)}{4p}}\nu^{-p\alpha}d\nu+\int_{\frac{1}{2}}^1(1-\nu)^{-\frac{3(p-1)}{4p}}\nu^{-p\alpha}d\nu\right].
\end{eqnarray*}
A simple computation gives
$$\int_0^{\frac{1}{2}}(1-\nu)^{-\frac{3(p-1)}{4p}}\nu^{-p\alpha}d\nu\leq\int_0^{\frac{1}{2}}\nu^{-\frac{3(p-1)}{4p}-p\alpha}d\nu=\frac{2^{\frac{3(p-1)}{4p}+p\alpha-1}}{1-\frac{3(p-1)}{4p}-p\alpha},$$
and
$$\int_{\frac{1}{2}}^1(1-\nu)^{-\frac{3(p-1)}{4p}}\nu^{-p\alpha}d\nu\leq\int_{\frac{1}{2}}^1(1-\nu)^{-\frac{3(p-1)}{4p}-p\alpha}d\nu=\frac{2^{\frac{3(p-1)}{4p}+p\alpha-1}}{1-\frac{3(p-1)}{4p}-p\alpha}.$$
As a consequence to these estimates
\begin{eqnarray*}
& & t^{\alpha}\int_0^t||S(t-s)(f(u_1+v)-f(u_2+v))||_{2p}ds\leq\\
& \leq & (4\pi)^{-\frac{3(p-1)}{4p}}C_f\max\{2^{p-2},1\}(2R_2+2T^{(p-1)\alpha}||v||_{\infty,2p}^{p-1})\times\nonumber\\
& &\times\frac{2^{\frac{3(p-1)}{4p}+p\alpha}}{1-\frac{3(p-1)}{4p}-p\alpha}T^{1+\alpha-p\alpha-\frac{3(p-1)}{4p}}||u_1-u_2||_{Y_T}\nonumber\\
& = & (4\pi)^{-\frac{3(p-1)}{4p}}C_f\max\{2^{p-2},1\}(2R_2+2T^{(p-1)\alpha}||v||_{\infty,2p}^{p-1})\times\nonumber\\
& & \times\frac{2^{\frac{3(p-1)}{4p}+p\alpha}}{1-\frac{3(p-1)}{4p}-p\alpha}T^{\frac{p+3}{8p}}||u_1-u_2||_{Y_T}.
\end{eqnarray*}
We put $C_1= (4\pi)^{-\frac{3(p-1)}{4p}}C_f\max\{2^{p-2},1\}(2R_2+2T^{(p-1)\alpha}||v||_{\infty,2p}^{p-1})\frac{2^{\frac{3(p-1)}{4p}+p\alpha}}{1-\frac{3(p-1)}{4p}-p\alpha}T^{\frac{p+3}{8p}}$, then $C_1\leq (4\pi)^{-\frac{3(p-1)}{4p}}C_f\max\{2^{p-2},1\}(2R_2+2T_*^{(p-1)\alpha}||v||_{\infty,2p}^{p-1})\frac{2^{\frac{3(p-1)}{4p}+p\alpha}}{1-\frac{3(p-1)}{4p}-p\alpha}T_*^{\frac{p+3}{8p}}\leq(4\pi)^{-\frac{3(p-1)}{4p}}C_f\max\{2^{p-2},1\}(2R_2+2||v||_{\infty,2p}^{p-1})\frac{2^{\frac{3(p-1)}{4p}+p\alpha}}{1-\frac{3(p-1)}{4p}-p\alpha}T_*^{\frac{p+3}{8p}}\leq\frac{1}{2}$. This implies that $\Phi$ is a contraction in the Banach space $Y_T$.
\\\h Choosing $u_2$ to be $0$ in this estimate, we obtain also that $||\Phi(u_1)-\Phi(0)||_{Y_T}\leq C_1||u_1||_{Y_T}\leq \frac{1}{2}R_2$. Moreover, the estimate
\begin{eqnarray}
 t^{\alpha}\int_0^t||S(t-s)(f(v))||_{2p}ds& \leq & \int_0^t(4\pi(t-s))^{-\frac{3(p-1)}{4p}}||f(v)||_{2}ds\nonumber\\
& \leq & \int_0^t(t-s)^{-\frac{3(p-1)}{4p}}||f(v)||_{2}ds\leq\frac{1}{2}R_2
\end{eqnarray}
implies that $||\Phi(u_1)||_{Y_T}\leq R_2$. Which means that $\Phi$ is a contraction from a complete metric space to itself and it admits a unique fixed-point $u$ in this set according to the Banach Theorem.
{\\\h\bf Step 2:} We prove that $u\in L^{\infty}({\mathcal O}\times(0,T))$. 
\\\h {\it Step 2.1:} We consider the case $v\in L^{\infty}(0,T,L^{2p}({\mathcal O}))$.
\\\h Firstly, we prove that $u\in L^{\infty}(0,T,L^{2p}({\mathcal O}))$. From $u(t,x)=\int_0^tS(t-s)f(u+v)(s,x)ds,$
we have the following inequality
\begin{eqnarray*}
||u(t,x)||_{2p}& \leq &\int_{0}^{t}||S(t-s)f(u+v)||_{2p}ds\leq\int_{0}^{t}(4\pi(t-s))^{-\frac{3(p-1)}{4p}}||f(u+v)||_{2}ds\\
& \leq &\int_{0}^{t}(4\pi(t-s))^{-\frac{3(p-1)}{4p}}(||f(u+v)-f(v)||_{2}+||f(v)||_{2})ds.
\end{eqnarray*}
The computations
\begin{eqnarray*}
||f(v)||_2 & \leq & ||f(0)||_2+||f'(\mu)v||_2\leq ||f(0)||_2+||f'(\mu)||_{\frac{2p}{p-1}}||v||_{2p}\\
& \leq &||f(0)||_2+||C_f|v|^{p-1}||_{\frac{2p}{p-1}}||v||_{2p}= ||f(0)||_2+C_f||v||_{2p}^{p}
\end{eqnarray*}
implies that $f(v)$ belongs to $L^{\infty}(0,T,L^{2}({\mathcal O}))$. Combining this with the estimate of $||u(t,x)||_{2p}$, we obtain 
\begin{eqnarray*}
& &||u(x,t)||_{2p}\leq\\
& \leq &\int_{0}^{t}(4\pi(t-s))^{-\frac{3(p-1)}{4p}}(||f'(v+\zeta)||_{\frac{2p}{p-1}}||u||_{2p}+||f(v)||_{\infty,2})ds\\
& \leq &\int_{0}^{t}(4\pi(t-s))^{-\frac{3(p-1)}{4p}}(||C_f|v+\zeta|^{p-1}||_{\frac{2p}{p-1}}||u||_{2p}+||f(v)||_{\infty,2})ds\\
& \leq &\int_{0}^{t}(4\pi(t-s))^{-\frac{3(p-1)}{4p}}(||C_f\max\{1,2^{p-2}\}(|v|^{p-1}+|u|^{p-1})||_{\frac{2p}{p-1}}||u||_{2p}\\
& &+||f(v)||_{\infty,2})ds\\
& \leq &\int_{0}^{t}(4\pi(t-s))^{-\frac{3(p-1)}{4p}}[C_f\max\{1,2^{p-2}\}(||v||_{2p}^{p-1}+||u||^{p-1}_{2p})||u||_{2p}\\
& &+||f(v)||_{\infty,2}]ds\\
& \leq &\left[\int_{0}^{t}(4\pi(t-s))^{-\frac{3(p-1)l_1}{4p}}ds\right]^{\frac{1}{l_1}}\times\\
& &\times\left[\int_0^t[C_f\max\{1,2^{p-2}\}(||v||_{\infty,2p}^{p-1}+||u||^{p-1}_{2p})||u||_{2p}+||f(v)||_{\infty,2}]^{l_2}ds\right]^{\frac{1}{l_2}},
\end{eqnarray*}
notice that $l_1$, $l_2$ $>0$, $\frac{1}{l_1}+\frac{1}{l_2}=1$ and $l_1<\frac{4p}{3(p-1)}$.
\\ The inequality
$$\int_{0}^{t}(t-s)^{-\frac{3(p-1)l_1}{4p}}ds=\left.-\frac{(t-s)^{1-\frac{3(p-1)l_1}{4p}}}{1-\frac{3(p-1)l_1}{4p}}\right|_0^t=\frac{t^{1-\frac{3(p-1)l_1}{4p}}}{1-\frac{3(p-1)l_1}{4p}}\leq\frac{T^{1-\frac{3(p-1)l_1}{4p}}}{1-\frac{3(p-1)l_1}{4p}}$$
leads to
\begin{eqnarray}\label{3e8}
& &||u||_{2p} \leq \left(\frac{T^{1-\frac{3(p-1)l_1}{4p}}}{1-\frac{3(p-1)l_1}{4p}}\right)^{\frac{1}{l_1}}(4\pi)^{-\frac{3(p-1)}{4p}}\times\\\nonumber
& &\times\left[\int_0^t[C_f\max\{1,2^{p-2}\}(||v||_{\infty,2p}^{p-1}+||u||^{p-1}_{2p})||u||_{2p}+||f(v)||_{\infty,2}]^{l_2}ds\right]^{\frac{1}{l_2}}.
\end{eqnarray}
Put $U(t)=||u(.,t)||_{2p}^{l_2}$, the inequality $(\ref{3e8})$ becomes
$$U(t)\leq \left(\frac{T^{1-\frac{3(p-1)l_1}{4p}}}{1-\frac{3(p-1)l_1}{4p}}\right)^{\frac{l_2}{l_1}}(4\pi)^{-\frac{3(p-1)l_2}{4p}}\times$$
$$\times\int_0^t[C_f\max\{1,2^{p-2}\}(||v||_{\infty,2p}^{p-1}+U(s)^{\frac{p-1}{l_2}})U(s)^{\frac{1}{l_2}}+||f(v)||_{\infty,2}]^{l_2}ds.$$
We denote $V(t)=\left(\frac{T^{1-\frac{3(p-1)l_1}{4p}}}{1-\frac{3(p-1)l_1}{4p}}\right)^{\frac{l_2}{l_1}}$ $(4\pi)^{-\frac{3(p-1)l_2}{4p}}$
 $\int_0^t[C_f\max\{1,2^{p-2}\}$ $(||v||_{\infty,2p}^{p-1}$ $+$ $U(s)^{\frac{p-1}{l_2}})$ $U(s)^{\frac{1}{l_2}}$ $+||f(v)||_{\infty,2}]^{l_2}ds$, then $V'(t)=\left(\frac{T^{1-\frac{3(p-1)l_1}{4p}}}{1-\frac{3(p-1)l_1}{4p}}\right)^{\frac{l_2}{l_1}}$ $(4\pi)^{-\frac{3(p-1)l_2}{4p}}$ $[C_f$ $\max\{1,2^{p-2}\}$ $(||v||_{\infty,2p}^{p-1}$ $+$ $U(t)^{\frac{p-1}{l_2}})$ $U(t)^{\frac{1}{l_2}}+||f(v)||_{\infty,2}]^{l_2}$. Since $U(t)\leq V(t)$, then $V'(t)\leq \left(\frac{T^{1-\frac{3(p-1)l_1}{4p}}}{1-\frac{3(p-1)l_1}{4p}}\right)^{\frac{l_2}{l_1}}$ $(4\pi)^{-\frac{3(p-1)l_2}{4p}}$ $[C_f$ $\max\{1,$ $2^{p-2}\}$ $(||v||_{\infty,2p}^{p-1}$ $+V(s)^{\frac{p-1}{l_2}})$ $V(s)^{\frac{1}{l_2}}$ $+$ $||f(v)||_{\infty,2}]^{l_2}$.
\\ Basing on the notations in $(\ref{2e1d})$, we put $G(r)=$
$$\int_0^r\frac{(4\pi)^{\frac{3(p-1)l_2}{4p}}d\zeta}{\left(\frac{T^{1-\frac{3(p-1)l_1}{4p}}}{1-\frac{3(p-1)l_1}{4p}}\right)^{\frac{l_2}{l_1}}[C_f\max\{1,2^{p-2}\}(||v||_{\infty,2p}^{p-1}+\zeta^{\frac{p-1}{l_2}})\zeta^{\frac{1}{l_2}}+||f(v)||_{\infty,2}]^{l_2}},$$
then $G(r)$ is an increasing function on $[0,+\infty)$ and $G'(V(t))\leq 1$. Hence $G(V(t))\leq t+G(V(0))=t+G(0)=t\leq T_*$ and $U(t)\leq V(t)\leq G^{-1}(T_*)$. Finally, $||u(.,t)||_{2p}\leq G^{-1}(T_*)^{\frac{1}{l_2}}$. We have proved that $u\in L^{\infty}(0,T,L^{2p}({\mathcal O}))$ by proving a {\it New Gronwall Inequality} on $U(t)$.
\\\h Secondly, we prove that $u\in L^{\infty}(0,T,L^{\infty}({\mathcal O}))$. Using Proposition 48.4 page 441, \cite{Quittner:2007:SPP}, we get the following estimate
\begin{eqnarray*}
& &||u(t)||_\infty\leq\\
& \leq & \int_0^t||S(t-s)f(u+v)(s)||_{\infty}ds\leq \int_0^t(4\pi(t-s))^{-\frac{3}{4}}||f(u+v)(s)||_{2}ds\\
& \leq &\int_0^t(4\pi(t-s))^{-\frac{3}{4}}(||f(u+v)(s)-f(v)||_2+||f(v)||_2)ds\\
& \leq &\int_0^t(4\pi(t-s))^{-\frac{3}{4}}[C_f\max\{1,2^{p-2}\}(||u||_{2p}^{p-1}+||v||^{p-1}_{\infty,2p})||u||_{2p}+||f(v)||_{\infty,2}]ds\\
& \leq &\int_0^t(4\pi(t-s))^{-\frac{3}{4}}[C_f\max\{1,2^{p-2}\}(G^{-1}(T_*)^{\frac{p-1}{l_2}}+||v||^{p-1}_{\infty,2p})G^{-1}(T_*)^{\frac{1}{l_2}}\\
& &+||f(v)||_{\infty,2}]ds\\
& \leq &\pi^{-\frac{3}{4}}(4T_*)^{\frac{1}{4}}[C_f\max\{1,2^{p-2}\}(G^{-1}(T_*)^{\frac{p-1}{l_2}}+||v||^{p-1}_{\infty,2p})G^{-1}(T_*)^{\frac{1}{l_2}}\\
& &+||f(v)||_{\infty,2}],
\end{eqnarray*}
which deduces $u\in L^{\infty}({\mathcal O}\times(0,T))$.
\\\h {\it Step 2.2:} For the case $v\in L^{\infty}(0,T,L^{\infty}({\mathcal O}))$, we use exactly the same argument, but with the function $G(r)=$
$$\int_0^r\frac{(4\pi)^{\frac{3(p-1)l_2}{4p}}d\zeta}{\left(\frac{T^{1-\frac{3(p-1)l_1}{4p}}}{1-\frac{3(p-1)l_1}{4p}}\right)^{\frac{l_2}{l_1}}\left[C_f\max\{1,2^{p-2}\}((Mm({\mathcal O})^{\frac{1}{2p}})^{p-1})+\zeta^{\frac{p-1}{l_2}})\zeta^{\frac{1}{l_2}}+{\max}_{|\zeta|\leq M}|f(\zeta)|m({\mathcal O})^{\frac{1}{2}}\right]^{l_2}}.$$
\\ Then we have $||u(t)||_\infty\leq$ $\pi^{-\frac{3}{4}}$ $(4T_*)^{\frac{1}{4}}$ $[C_f$ $\max\{1,$ $2^{p-2}\}$ $(G^{-1}$ $(T_*)^{\frac{p-1}{l_2}}$ $+(Mm$ $({\mathcal O})^{\frac{1}{2p}})^{p-1})$ $G^{-1}(T_*)^{\frac{p-1}{l_2}}$ $+{\max}_{|\zeta|\leq M}$ $|f(\zeta)|$ $m({\mathcal O})^{\frac{1}{2}}]$ and $u\in L^{\infty}({\mathcal O}\times(0,T))$. 
{\\\h\bf Step 3:} We prove that $u\in C([0,T],L^2({\mathcal O}))$ and $u$ is also a solution of $(\ref{2e1})$. The proof in this step works well for both cases  $v\in L^{\infty}(0,T,L^{2p}({\mathcal O}))$ and $v\in L^{\infty}(0,T,L^{\infty}({\mathcal O}))$.
\\\h For $\epsilon$ positive, 
\begin{eqnarray*}
& &u(.,t+\epsilon)-u(.,t)\\
& = &\int_0^{t+\epsilon}S(t+\epsilon-s)f(u+v)(s)ds-\int_0^{t}S(t-s)f(u+v)(s)ds\\
& = &\int_0^t[S(t+\epsilon-s)-S(t-s)]f(u+v)(s)ds+\int_t^{t+\epsilon}S(t+\epsilon-s)f(u+v)(s)ds,
\end{eqnarray*}
which implies the following inequality
\begin{eqnarray*}
||u(.,t+\epsilon)-u(.,t)||_2& \leq &\int_0^t||[S(t+\epsilon-s)-S(t-s)]f(u+v)(s)||_2ds+\\
& &+\int_t^{t+\epsilon}||S(t+\epsilon-s)f(u+v)(s)||_2ds.
\end{eqnarray*}
\h We firstly estimate the second term on the RHS of the previous inequality. Due to Proposition 48.4 page 441, \cite{Quittner:2007:SPP}, 
$$\int_t^{t+\epsilon}||S(t+\epsilon-s)f(u+v)(s)||_2ds\leq \int_t^{t+\epsilon}||f(u+v)(s)||_2ds.$$
Since $\int_t^{t+\epsilon}||f(u+v)(s)||_2ds$ tends to $0$ as $\epsilon$ tends to $0$, $\int_t^{t+\epsilon}||S(t+\epsilon-s)f(u+v)(s)||_2ds$ tends to $ 0$ as $\epsilon$ tends to $0$.
\\\h Now, we estimate the first term. Noticing that $||[S(t+\epsilon-s)-S(t-s)]f(u+v)(s)||_2$ tends to $0$ as $\epsilon$ tends to $0$ for all $s$ in $[0,t]$, and
\begin{eqnarray*}
& &||[S(t+\epsilon-s)-S(t-s)]f(u+v)(s)||_2\\
& \leq &||S(t+\epsilon-s)f(u+v)(s)||_2+||S(t-s)f(u+v)(s)||_2\leq 2||f(u+v)(s)||_2,
\end{eqnarray*}
where $f(u+v)$ can be proved to belong to $L^{\infty}(0,T,L^2({\mathcal O}))$ by using the same argument as before; due to the Lebesgue Dominated Convergence Theorem, we deduce that
$$\lim_{\epsilon\to 0}\int_0^t||[S(t+\epsilon-s)-S(t-s)]f(u+v)(s)||_2ds=0.$$
\\\h Therefore $\lim_{\epsilon\to 0}||u(.,t+\epsilon)-v(.,t)||_2=0$, or $u\in C([0,T],L^2({\mathcal O}))$.
\\\h Since $v\in C([0,T],L^2({\mathcal O}))$, $u+v\in C([0,T],L^2({\mathcal O}))$. Which implies $f(u+v)\in C([0,T],L^2({\mathcal O}))$; and by Lemma $\ref{3l2}$, $u$ is also a solution of $(\ref{2e1})$.
\section{Proof of The Existence Results for Semilinear Heat Equation in a Cylindrical Bounded Domain - Theorem $\ref{2t2}$} First, we consider the following equation
\begin{equation}
\left \{ \begin{array}{ll}\partial_tv-\Delta v=0&\mbox{ in }\Omega\times(0,T_*),\vspace{.1in}\\ 
v=g&\mbox{ on }\partial\Omega,\vspace{.1in}\\ 
v(.,0)=u_0(.)&\mbox{ on }\bar\Omega.\end{array}\right. 
\label{4e1}
\end{equation}
\\ We recall the following definition:
\begin{definition}(Definition c page 69, \cite{Friedman:1964:PDE}) Let $\omega$ be a domain in $\mathbb{R}^N$, we say that $\partial\omega\times(0,T]$ has the outside strong sphere property if for every $Q=(x_0,t_0)$ in $\partial\omega\times(0,T]$, there exists a ball $K$ with center $(\bar x,\bar t)$ such that $K\cap(\bar\omega\times[0,T])=\{Q\}$ and $|\bar x-x|\geq\mu(Q)>0$, for every $(x,t)$ in $\bar\omega\times[0,T]$, $|t-t_0|<\epsilon$ where $\epsilon$ is a constant independent of $Q$. 
\label{4d1}
\end{definition}
We will prove that, in fact $\partial\Omega\times(0,T]$ has the outside strong sphere property. 
\\ Let $Q=(x_0,t_0)$ be on $\partial\Omega\times(0,T]$. Since $D$ is smooth, there exists $\bar x$ and a constant $r>0$ such that $K=\bar B_{\mathbb{R}^N}(\bar x,r)\cap\bar\Omega=\{x_0\}$ and $||\bar x-x_0||_{\mathbb{R}^N}=r$. 
\\ We choose $\bar t=t_0$ and let $(x',t')$ be in $\overline{B_{\mathbb{R}^{N+1}}((\bar x,\bar t),r)}\cap(\bar\Omega\times[0,T])$, then $||x'-\bar x||_{\mathbb{R}^3}+|t'-\bar t|\leq r$. Since $||x'-\bar x||_{\mathbb{R}^N}\leq r$; $x'\in \overline{B_{\mathbb{R}^N}(\bar x,r)}\cap\bar\Omega$ and $x'=x_0$. We infer that $t'=t_0$ and $\overline{B_{\mathbb{R}^{N+1}}((\bar x,\bar t),r)}\cap(\bar\Omega\times[0,T])=\{(x_0,t_0)\}.$ For $(x,t)\in\bar\Omega\times[0,T]$, $||x-\bar x||_{\mathbb{R}^N}\geq r>0.$ Consequently, $\partial\Omega\times(0,T]$ has the outside strong sphere property. 
\\ By Theorem 8 page 69, \cite{Friedman:1964:PDE}, we can conclude that $\partial\Omega\times(0,T]$ has local barriers. Due to Theorem 5 page 123, \cite{Friedman:1964:PDE}, there exists a unique solution $v$ to the problem $(\ref{4e1})$, and $v\in C^{1,2}((0,T]\times\Omega)\cap C(\bar\Omega\times[0,T])$. According to Theorem 10 page 72, \cite{Friedman:1964:PDE}, $v\in C^{\infty}(\Omega\times(0,T))$.
\\\h Due to the Week Maximum Principle page 368, \cite{Evans:1998:PDE}: $\max_{\bar\Omega\times[0,T]}v$ $=$ $\max_{\partial(\Omega\times(0,T))}$ $v,$ and $\min_{\bar\Omega\times[0,T]}v=\min_{\partial(\Omega\times(0,T))}v.$ Which gives $v<M$ and $v>-M$ on $\bar\Omega\times[0,T]$ since $M>\max\{||u_0||_{\infty},$ $||g||_{\infty}\}$. Then $||v||_{L^{\infty}(\Omega\times(0,T))}<M.$
\\\h Subtracting $(\ref{4e1})$ and $(\ref{2e2})$, and put $w=u-v$, we have the following equations
\begin{equation}
\left \{ \begin{array}{ll}\partial_tw-\Delta w=f(w+v)&\mbox{ in }\Omega\times(0,T),\vspace{.1in}\\ 
w(.,.)=0&\mbox{ on }\partial\Omega,\vspace{.1in}\\ 
w(.,0)=0&\mbox{ on }\bar\Omega.\end{array}\right. 
\label{4e2}
\end{equation}
According to Theorem $\ref{2t1}$, $(\ref{4e2})$ has a solution $w$ satisfying $w\in L^{\infty}(\Omega\times(0,T))\cap C([0,T],L^2(\Omega))\cap L^{2}(0,T,H_0^1(\Omega))$ and $\partial_t w\in L^{2}(0,T,L^2(\Omega))$. Moreover, $||w||_{\infty,\infty}$ $\leq$ $\pi^{-\frac{3}{4}}$ $(4T_*)^{\frac{1}{4}}$ $[C_f$ $\max\{1,2^{p-2}\}$ $(G^{-1}(T_*)^{\frac{p-1}{l_2}}$ $+((Mm(\Omega)^{\frac{1}{2p}})^{p-1})$ $G^{-1}(T_*)^{\frac{p-1}{l_2}}$ $+{\max}_{|\zeta|\leq M}$ $|f(\zeta)|m(\Omega)^{\frac{1}{2}}],$ where $G(r)=$
$$\int_0^r\frac{(4\pi)^{\frac{3(p-1)l_2}{4p}}d\zeta}{\left(\frac{T^{1-\frac{3(p-1)l_1}{4p}}}{1-\frac{3(p-1)l_1}{4p}}\right)^{\frac{l_2}{l_1}}\left[C_f\max\{1,2^{p-2}\}((Mm(\Omega)^{\frac{1}{2p}})^{p-1}+\zeta^{\frac{p-1}{l_2}})\zeta^{\frac{1}{l_2}}+{\max}_{|\zeta|\leq M}|f(\zeta)|m(\Omega)^{\frac{1}{2}}\right]^{l_2}}.$$
\h Hence $(\ref{2e2})$ has a solution $u$ satisfying $u\in$ $ L^{\infty}(\Omega\times(0,T))$ $\cap$ $ C([0,T],L^2(\Omega))$ $\cap$ $ L^{2}(0,T,H_0^1(\Omega))$ and $\partial_t u\in L^{2}(0,T,L^2(\Omega))$. Moreover, $||u||_{\infty,\infty}$ $\leq$ $M_*:=$$M+$ $\pi^{-\frac{3}{4}}$ $(4T_*)^{\frac{1}{4}}$ $[C_f$ $\max\{1,2^{p-2}\}$ $(G^{-1}(T_*)^{\frac{p-1}{l_2}}$ $+((Mm(\Omega)^{\frac{1}{2p}})^{p-1})$ $G^{-1}(T_*)^{\frac{p-1}{l_2}}$ $+{\max}_{|\zeta|\leq M}$ $|f(\zeta)|m(\Omega)^{\frac{1}{2}}].$ 
\\\h Using the results in section VI.8 \cite{Lieberman:1996:SOP} for the equation $\partial_tu-\Delta u=g$ in $\Omega\times(0,T)$, where $g=f(u+v)$, $g\in L^{\infty}(\Omega\times(0,T))$, we can conclude that $u$ is continuous on $\overline{\Omega\times(0,T)}$.
\section{Proof of the Well-Posedness and Convergence Properties of the Algorithms - Theorems $\ref{2t3}$ and $\ref{2t4}$}
\subsection{The Well-Posedness of the Algorithm - Proof of Theorem $\ref{2t3}$} We consider the equation $(\ref{2e5})$
\begin{equation}
\left \{ \begin{array}{ll}\partial_t\phi_M-\Delta\phi_M=f(\phi_M)&\mbox{ in }\Omega\times(0,T_0),\vspace{.1in}\\ 
\phi_M(.,.)=M&\mbox{ on }\partial\Omega\times[0,T_0],\vspace{.1in}\\ 
\phi_M(.,0)=M&\mbox{ on }\bar\Omega.\end{array}\right. 
\label{5e2}
\end{equation}
Since $f(x)$ is positive for all $x$ in $\mathbb{R}$, then $\partial_t\phi_M-\Delta\phi_M\geq 0$ and $\partial_t(\phi_M-M)-\Delta(\phi_M-M)\geq0$. Since $\phi_M=M$ on $(\partial\Omega\times[0,T_0])\cup(\bar\Omega\cap\{0\})$, then $\phi_M\geq M$ on $\overline{\Omega\times[0,T_0]}$ according to the Maximum Principle.
\\\h We will prove that the algorithm is well posed for $T<T^*$ by recursion.
\\\h For $k=1$, we consider the problem on the $j$-th domain
\begin{equation}
\left \{ \begin{array}{ll}
\partial_tu_{j}^{1}-\Delta u_{j}^{1}=f(u_{j}^{1})&\mbox{ in }\Omega_j\times(0,T),\vspace{.1in}\\
u_{j}^{1}(\cdot,a_j,\cdot)=u_{j-1}^{0}(\cdot,a_j,\cdot) &\mbox{ in }D\times(0,T),\vspace{.1in}\\
u_{j}^{1}(\cdot,b_j,\cdot)=u_{j+1}^{0}(\cdot,b_j,\cdot) &\mbox{ in }D\times(0,T),\vspace{.1in}\\
u_{j}^{1}(.,.)=g(.,.)&\mbox{ on }\partial D\times[a_j,b_j]\times[0,T],\vspace{.1in}\\ 
u_{j}^{1}(.,0)=u_0(.)&\mbox{ on }\overline{\Omega_j}.\end{array}\right. 
\label{5e1}
\end{equation}
According to Theorem $\ref{2t2}$, the equation $(\ref{5e1})$ has a solution $u_{j}^{1}$ satisfying $u_{j}^{1}\in L^{\infty}(\Omega_j\times(0,T))\cap C([0,T],L^2(\Omega_j))\cap L^{2}(0,T,H_0^1(\Omega_j))$, $\partial_t u_{j}^{1}\in L^{2}(0,T,L^2(\Omega_j))$and $u_{j}^{1}$ is continuous on $\overline{\Omega_j\times(0,T)}$. The fact that $u_{j}^{1}\in C^{1,2}((0,T)\times\Omega))$ then follows.
\\ Since $\partial_t(\phi_M-u_{j}^{1})-\Delta(\phi_M-u_{j}^{1})=f(\phi_M)-f(u_{j}^{1})$; $\partial_t(\phi_M-u_{j}^{1})-\Delta(\phi_M-u_{j}^{1})-(\phi_M-u_{j}^{1})c(x,t)=0$, where
\begin{eqnarray*}
c(x,t)=\left \{ \begin{array}{ll}\frac{f(\phi_M)-f(u_{j}^{1})}{\phi_M-u_{j}^{1}} &\mbox{ if }\phi_M\ne u_{j}^{1} ,\vspace{.1in}\\ 
0&\mbox{ otherwise }.\end{array}\right. 
\end{eqnarray*}
We notice that $c\in C(\overline{\Omega_j\times(0,T)})$ and $u_{j}^{1}\leq M\leq\phi_M$ on $(\partial\Omega_j\times[0,T])\cup(\bar\Omega_j\times\{0\})$. Then it follows by the Maximum Principle that $\phi_M\geq u_{j}^{1}$ on $\overline{\Omega_j\times(0,T)}$.
\\ Since $\partial_t(u_{j}^{1}+M)-\Delta(u_{j}^{1}+M)-(f(u_{j}^{1})-f(-M))=f(-M)\geq 0$ and $u_{j}^{1}\geq -M$ on $(\partial\Omega_j\times[0,T])\cup(\bar\Omega_j\times\{0\})$; the Maximum Principle then implies that $u_{j}^{1}\geq -M$ on $\overline{\Omega_j\times(0,T)}$.
\\ Consequently, $\phi_M\geq u_{j}^{1}\geq -M$ on $\overline{\Omega_j\times(0,T)}$, $\forall j=\overline{1,I}$.
\\\h Suppose that the algorithm is well posed up to step $k=m$, and moreover $\phi_M\geq u_{j}^{i}\geq -M$ on $\overline{\Omega_j\times(0,T)}$ $\forall j=\overline{1,I}, i=\overline{1,m}$, using the same argument as above we can see that the algorithm is well posed for $k=m+1$ and $\phi_M\geq u_{j}^{m+1}\geq -M$ on $\overline{\Omega_j\times(0,T)}$, $\forall j=\overline{1,I}$.
\\\h Therefore, the algorithm is well posed for all $k\in\mathbb{N}$ and we have also $\phi_M\geq u_{j}^k\geq -M$ on $\overline{\Omega_j\times(0,T)}$, $\forall j=\overline{1,I}$.
\subsection{The Convergence of the Algorithm - Proof of Theorem $\ref{2t4}$} We devide the proof into several steps.
{\\\h\bf Step 1:} The Exponential Decay Error Estimates.
\\\h On the $j$-th domain, at the $k$-th iteration, we consider the equation
\begin{equation}
\left \{ \begin{array}{ll}\partial_te_{j}^k-\Delta e_{j}^k=f(u_j^k)-f(u)&\mbox{ in }\Omega_j\times(0,T),\vspace{.1in}\\ 
e_{j}^k(\cdot,a_j,\cdot)=e_{j-1}^{k-1}(\cdot,a_j,\cdot)&\mbox{ in }D\times(0,T),\vspace{.1in}\\
e_{j}^k(\cdot,b_j,\cdot)=e_{j+1}^{k-1}(\cdot,b_j,\cdot)&\mbox{ in }D\times(0,T).\vspace{.1in}\\
e_{j}^k(.,.)=0&\mbox{ on }\partial D\times[a_j,b_j]\times[0,T],\vspace{.1in}\\ 
e_{j}^k(.,0)=0&\mbox{ on }\bar\Omega_j.\end{array}\right. 
\label{5e3}
\end{equation}
\\\h Since the algorithm is well posed, according to Definition $\ref{2d1}$, the set $\{u_j^k|j=\overline{1,I},k\in\mathbb{N}\}$ is bounded by a constant $C_1$. Which implies that $\{e_{j}^k|j=\overline{1,I},k\in\mathbb{N}\}$ is bounded by a constant $C_2$, where $C_2=C_1+||u||_{C(\overline{\Omega\times(0,T)})}$.
\\\h We consider the Controlling Function $\Phi(x,t)=P(e_{j}^k)\exp(\beta(z-a_{j})-\gamma t),$  and ${\cal L}(\Phi)=\partial_t\Phi-\Delta\Phi+2\beta\partial_z\Phi$, where $\beta$ and $\gamma$ are constants to be chosen later in the next steps. We now prove that if we choose $(\beta,\gamma)$ such that $\gamma-\beta^2$ is large enough, ${\cal L}(\Phi)$ is negative.
\\ A simple computation gives 
$${\cal L}(\Phi)=(\partial_te_{j}^k-\Delta e_{j}^k)P'(e_{j}^k)\exp(\beta(z-a_{j})-\gamma t)+$$
$$+(\beta^2-\gamma)P(e_{j}^k)\exp(\beta(z-a_{j})-\gamma t)-(\nabla e_{j}^k)^2P''(e_{j}^k)\exp(\beta(z-a_{j})-\gamma t).$$
By $(\ref{5e3})$, the previous computation leads to
$${\cal L}(\Phi)=(f(u_j^k)-f(u))P'(e_{j}^k)\exp(\beta(z-a_{j})-\gamma t)+$$
$$+(\beta^2-\gamma)P(e_j^k)\exp(\beta(z-a_{j})-\gamma t)-(\nabla e_j^k)^2P''(e_j^k)\exp(\beta(z-a_{j})-\gamma t),$$
which implies
$${\cal L}(\Phi)\leq[(f(u_j^k)-f(u))P'(e_j^k)+(\beta^2-\gamma)P(e_j^k)]\exp(\beta(z-a_{j})-\gamma t).$$
We consider the term $(f(u_j^k)-f(u))P'(e_j^k)+(\beta^2-\gamma)P(e_j^k)$. By the Mean Value Theorem, $f(u_j^k(x,t))-f(u(x,t))=e_j^k(x,t)f'(\zeta(x,t))$, where $\zeta(x,t)$ is a number lying between $u_j^k(x,t)$ and $u(x,t)$. We observe that $|\zeta(x,t)|$ $\leq$ $ |u_j^k(x,t)|$ $+|u(x,t)|\leq C_1+||u||_{C(\overline{\Omega\times(0,T)})}=C_2$, which implies $|f'(\zeta(x,t))|$ is bounded on $\overline{\Omega\times(0,T)}$. 
\\ We fix a pair $(x,t)$ in $\Omega\times(0,T)$: 
\begin{itemize}
\item If $P(e_j^k(x,t))=0$, then $e_j^k=0$ and $P'(e_j^k)=0$. That means $(f(u_j^k)-f(u))P'(e_j^k)+(\beta^2-\gamma)P(e_j^k)=0$. 
\item If $P(e_j^k(x,t))\ne0$, then 
\begin{eqnarray*}
\frac{(f(u_j^k)-f(u))P'(e_j^k)+(\beta^2-\gamma)P(e_j^k)}{P(e_j^k)}=\frac{e_j^k(x,t)P'(e_j^k)}{P(e_j^k)}f'(\zeta(x,t))+(\beta^2-\gamma).
\end{eqnarray*}
 Since $e_j^k$ is bounded by $C_2$, $\left|\frac{e_j^k(x,t)P'(e_j^k)}{P(e_j^k)}\right|$ is bounded by $K(C_2)$; which means $\frac{e_j^k(x,t)P'(e_j^k)}{P(e_j^k)}$ $f'(\zeta(x,t))$ is bounded on $\overline{\Omega\times(0,T)}$ as $|f'(\zeta(x,t))|$ is bounded. Since $P(e_j^k)$ is positive, then if we choose $(\beta,\gamma)$ such that $-\beta^2+\gamma$ is large enough, we can have that $(f(u_j^k)-f(u))P'(e_j^k)+(\beta^2-\gamma)P(e_j^k)$ is negative.
\end{itemize}
Consequently, if $\gamma-\beta^2$ is large enough, $(f(u_{j,k})-f(u))P'(e_j^k)+(\beta^2-\gamma)P(e_j^k)$ is negative or ${\cal L}(\Phi)$ is negative.
\\\h When ${\cal L}(\Phi)\leq 0$, according to the Maximum Principle, the function $\Phi(x,t)$ can only attain its maximum values on the boundary $\partial\Omega\times[0,T]$ or $\bar\Omega\times\{0\}$. We have that $\Phi\geq 0$, and moreover $\Phi=0$ on $\bar\Omega_j\times\{0\}$ and on $\partial D\times[a_j,b_j]\times[0,T]$. Thus:
\begin{itemize}
\item If the maximum value(s) of $\Phi$ can be achived on $\bar\Omega_j\times\{0\}$ and on $\partial D\times[a_j,b_j]\times[0,T]$, then $\Phi=0$ on $\overline{\Omega\times(0,T)}$; which means that $e_j^k=0$ on $\overline{\Omega\times(0,T)}$. 
\item If $e_j^k\ne0$, then the maximum value(s) of $\Phi$ can be achived only on $\overline{D}\times{\{a_j\}}\times[0,T]$ or on $\overline{D}\times{\{b_j\}}\times[0,T]$. 
\end{itemize}
\h For $x$ in $\mathbb{R}^N$, denote that $x=(X,z)$ where $X\in\mathbb{R}^{N-1}$ and $z\in\mathbb{R}$, we now have the following exponential decay estimates for the errors $e_j^k$.
{\\\h\it Case 1}: $j=1$.
\\\h We have that $e_{1}^{k}=0$ on $\overline{D}\times{\{a_1\}}\times[0,T]$. Since $e_1^k\ne0$, the maximum value(s) of $\Phi$ can only be achived on $\overline{D}\times{\{b_1\}}\times[0,T]$ and for $(X,z,t)\in\overline{\Omega_1\times(0,T)}$
\begin{eqnarray}
& &P(e_{1}^{k})(X,z,t)\exp(\beta(z-a_1)-\gamma t)\leq\nonumber\\ 
& \leq & \max_{\bar D\times[0,T]}\{P(e_{1}^{k})(X,b_1,t)\exp(\beta(b_1-a_1)-\gamma t)\}\nonumber\\
& = &\max_{\bar D\times[0,T]}\{P(e_{2}^{k-1})(X,b_1,t)\exp(\beta S_{1} -\gamma t)\}.
\label{5e4o1}
\end{eqnarray}
{\\\h\it Case 2}: $j=I$.
\\\h We have that $e_{I}^{k}=0$ on $\overline{D}\times{\{b_I\}}\times[0,T]$. Since $e_I^k\ne0$, the maximum value(s) of $\Phi$ can only be achived on $\overline{D}\times{\{a_I\}}\times[0,T]$ and for $(X,z,t)\in\overline{\Omega_I\times(0,T)}$
\begin{eqnarray}
& & P(e_{I}^{k})(X,z,t)\exp(\beta(z-a_{I})-\gamma t)\leq \nonumber\\
& \leq & \max_{\bar D\times[0,T]}\{P(e_{I}^{k})(X,a_{I},t)\exp(-\gamma t)\}\nonumber\\
& = &\max_{\bar D\times[0,T]}\{P(e_{I-1}^{k-1})(X,a_{I},t)\exp(-\gamma t)\}.
\label{5e4o2}
\end{eqnarray}
{\\\h\it Case 3}: $1<j<I$.
\\\h The maximum value(s) of $\Phi$ can be achived on both $\overline{D}\times{\{a_j\}}\times[0,T]$ and $\overline{D}\times{\{b_j\}}\times[0,T]$ and for $(X,z,t)\in\overline{\Omega_j\times(0,T)}$
\begin{eqnarray}
& &P(e_j^k)(X,z,t)\exp(\beta(z-a_{j})-\gamma t) \leq\nonumber\\
& \leq &\max\{\max_{\bar D\times[0,T]}\{P(e_j^k)(X,b_{j},t)\exp(\beta(b_{j}-a_{j})-\gamma t)\},\nonumber\\
& &\max_{\bar D\times[0,T]}\{P(e_j^k)(X,a_{j},t)\exp(-\gamma t)\}\}\nonumber\\
& = &\max\{\max_{\bar D\times[0,T]}\{P(e_{j+1}^{k-1})(X,b_{j},t)\exp(\beta(b_{j}-a_{j})-\gamma t)\},\nonumber\\
& &\max_{\bar D\times[0,T]}\{P(e_{j-1}^{k-1})(X,a_{j},t)\exp(-\gamma t)\}\}.
\label{5e4o3}
\end{eqnarray}
{\h\bf Step 2:} Proof of convergence.
{\\\h\bf Step 2.1:} Estimate of the right boudaries of the sub-domains.
\\\h Consider the $I$-th domain, at the $k$-th step, $(\ref{5e4o2})$ infers
$$P(e_I^k(X,z,t))\exp(\beta(z-a_{I})-\gamma t)\leq \max_{\bar D\times[0,T]}\{P(e_I^k(X,a_{I},t))\exp(-\gamma t)\}.$$
Replace $z$ by $b_{I-1}$, we get
$$P(e_I^k(X,b_{I-1},t))\exp(\beta(b_{I-1}-a_{I})-\gamma t)\leq \max_{\bar D\times[0,T]}\{P(e_I^k(X,a_{I},t))\exp(-\gamma t)\}.$$
Since $e_I^k(X,b_{I-1},t)=e_{I-1}^{k+1}(X,b_{I-1},t)$, 
$$P(e_{I-1}^{k+1}(X,b_{I-1},t))\exp(\beta(b_{I-1}-a_{I})-\gamma t)\leq \max_{\bar D\times[0,T]}\{P(e_I^k(X,a_{I},t))\exp(-\gamma t)\}.$$
Let $\beta$ in this case be $\beta_1$ to be $\sqrt\frac{\gamma}{2}$, then when we choose $\gamma$ large, $\gamma-\beta^2$ is large. The inequality becomes
$$P(e_{I-1}^{k+1}(X,b_{I-1},t))\exp(-\gamma t)\leq\exp(-\beta_1S_{I-1})\max_{\bar D\times[0,T]}\{P(e_I^k(X,a_{I},t))\exp(-\gamma t)\}.$$
We deduce
\begin{equation}
P(e_{I-1}^{k+1}(X,b_{I-1},t))\exp(-\gamma t)\leq \exp(-\beta_1S_{I-1})E_k.
\label{5e4}
\end{equation}
\h Moreover, on the $(I-1)$-th domain, at the $(k+1)$-th step, $(\ref{5e4o3})$ leads to
\begin{eqnarray*}
& &P(e_{I-1}^{k+1})(X,z,t)\exp(\beta(z-a_{I-1})-\gamma t)\leq\\
& & \max\{\max_{\bar D\times[0,T]}\{P(e_{I-1}^{k+1})(X,b_{I-1},t)\exp(\beta(b_{I-1}-a_{I-1})-\gamma t)\},\\
& &\max_{\bar D\times[0,T]}\{P(e_{I-1}^{k+1})(X,a_{I-1},t)\exp(-\gamma t)\}\}.
\end{eqnarray*}
Since $e_{I-1}^{k+1}(X,b_{I-2},t)=e_{I-2}^{k+2}(X,b_{I-2},t)$, 
\begin{eqnarray*}
& &P(e_{I-2}^{k+2})(X,b_{I-2},t)\exp(\beta(b_{I-2}-a_{I-1})-\gamma t)\leq\\
& \leq &\max\{\max_{\bar D\times[0,T]}\{P(e_{I-1}^{k+1})(X,b_{I-1},t)\exp(\beta(b_{I-1}-a_{I-1})-\gamma t)\},\\
& &\max_{\bar D\times[0,T]}\{P(e_{I-1}^{k+1})(X,a_{I-1},t)\exp(-\gamma t)\}\}\\
& = &\max\{\max_{\bar D\times[0,T]}\{P(e_{I-1}^{k+1})(X,b_{I-1},t)\exp(\beta L_{I-1}-\gamma t)\},\\
& &\max_{\bar D\times[0,T]}\{P(e_{I-1}^{k+1})(X,a_{I-1},t)\exp(-\gamma t)\}\}.
\end{eqnarray*}
Hence
\begin{eqnarray*}
& &P(e_{I-2}^{k+2})(X,b_{I-2},t)\exp(\beta S_{I-2}-\gamma t)\leq\max\{P(e_{I-1}^{k+1})(X,b_{I-1},t)\times\\
& &\times\exp(\beta L_{I-1}-\gamma t),P(e_{I-1}^{k+1})(X,a_{I-1},t)\exp(-\gamma t)\}.
\end{eqnarray*}
Since $$P(e_{I-1}^{k+1})(X,a_{I-1},t)\exp(-\gamma t)\leq E_{k+1},$$ $(\ref{5e4})$ implies  
$$P(e_{I-2}^{k+2})(X,b_{I-2},t)\exp(\beta S_{I-2}-\gamma t)\leq\max\{E_{k}\exp(\beta L_{I-1}-\beta_1S_{I-1}),E_{k+1}\}.$$
Thus
$$P(e_{I-2}^{k+2})(X,b_{I-2},t)\exp(-\gamma t)\leq\max\{E_{k}\exp(\beta (L_{I-1}-S_{I-2})-\beta_1S_{I-1}),E_{k+1}\exp(-\beta S_{I-2})\}.$$
We choose $\beta=\beta_2=\beta_1\frac{S_{I-1}}{L_{I-1}}$ such that $\beta_2(-L_{I-1}+S_{I-2})+\beta_1S_{I-1}=\beta_2S_{I-2}$; then
\begin{equation}
P(e_{I-2}^{k+2})(X,b_{I-2},t)\exp(-\gamma t)\leq\max\{E_{k},E_{k+1}\}\exp(-\beta_2S_{I-2}).
\label{5e5}
\end{equation}
\h Using the same techniques as the ones that we use to achive $(\ref{5e4})$ and $(\ref{5e5})$, we can prove that
\begin{equation}
P(e_{I-j}^{k+j})(X,b_{I-j},t)\exp(-\gamma t)\leq\max\{E_{k},\dots,E_{k+j-1}\}\exp(-\beta_jS_{I-j}), 
\label{5e6}
\end{equation}
where $\beta_j=\beta_1\frac{S_{I-1}}{L_{I-1}}\dots\frac{S_{I-j+1}}{L_{I-j+1}}$, $j=\{2,\dots,I-1\}$.
\\ {\h\bf Step 2.2:} Estimate of the left boundaries of the sub-domains
\\\h Consider the $1$-th domain, at the $k$-th step, $(\ref{5e4o1})$ infers
$$P(e_1^k(X,z,t))\exp(\beta z-\gamma t)\leq \max_{\bar D\times[0,T]}\{P(e_1^k(X,b_{1},t))\exp(\beta b_{1}-\gamma t)\}.$$
Replace $z$ by $a_{2}$, we get
$$P(e_1^k(X,a_{2},t))\exp(-\gamma t)\leq \max_{\bar D\times[0,T]}\{P(e_1^k(X,b_{1},t))\exp(-\gamma t)\}\exp(\beta(-a_{2}+b_{1})).$$
Since $e_1^k(X,a_{2},t)=e_2^{k+1}(X,a_{2},t)$, 
$$P(e_2^{k+1}(X,a_{2},t))\exp(-\gamma t)\leq \max_{\bar D\times[0,T]}\{P(e_1^k(X,b_{1},t))\exp(-\gamma t)\}\exp(\beta S_{1}).$$
Let $\beta=-\beta_1'=-\sqrt\frac{\gamma}{2}$, then when we choose $\gamma$ large, $\gamma-\beta^2$ is large. The inequality becomes
$$P(e_2^{k+1}(X,a_{2},t))\exp(-\gamma t)\leq\exp(-\beta'_1S_1)\max_{\bar D\times[0,T]}\{P(e_1^k(X,b_{1},t))\exp(-\gamma t)\}.$$
We deduce
\begin{equation}
P(e_{2}^{k+1}(X,a_{2},t))\exp(-\gamma t)\leq \exp(-\beta'_1S_1)E_k.
\label{5e7}
\end{equation}
\h Moreover, on the $2$-th domain, at the $(k+1)$-th step, $(\ref{5e4o3})$ leads to
\begin{eqnarray*}
& &P(e_{2}^{k+1})(X,z,t)\exp(\beta(z-a_{2})-\gamma t)\leq\\
& & \max\{\max_{\bar D\times[0,T]}\{P(e_{2}^{k+1})(X,b_{2},t)\exp(\beta(b_{2}-a_{2})-\gamma t)\},\\
& &\max_{\bar D\times[0,T]}\{P(e_{2}^{k+1})(X,a_{2},t)\exp(-\gamma t)\}\}.
\end{eqnarray*}

Since $e_{2}^{k+1}(X,a_{3},t)=e_{3}^{k+2}(X,a_{3},t)$, 
\begin{eqnarray*}
& &P(e_{3}^{k+2})(X,a_{3},t)\exp(\beta(a_{3}-a_{2})-\gamma t)\leq\\
& \leq &\max\{\max_{\bar D\times[0,T]}\{P(e_{2}^{k+1})(X,b_{2},t)\exp(\beta(b_{2}-a_{2})-\gamma t)\},\\
& &\max_{\bar D\times[0,T]}\{P(e_{2}^{k+1})(X,a_{2},t)\exp(-\gamma t)\}\}\\
& = &\max\{\max_{\bar D\times[0,T]}\{P(e_{2}^{k+1})(X,b_{2},t)\exp(\beta L_{2}-\gamma t)\},\\
& &\max_{\bar D\times[0,T]}\{P(e_{2}^{k+1})(X,a_{2},t)\exp(-\gamma t)\}\}.
\end{eqnarray*}
Hence
$$P(e_{3}^{k+2})(X,a_{3},t)\exp(\beta(L_2-S_{2})-\gamma t)\leq\max\{P(e_{2}^{k+1})(X,b_{2},t)\times$$
$$\times\exp(\beta L_{2}-\gamma t),P(e_{2}^{k+1})(X,a_{2},t)\exp(-\gamma t)\}.$$
Since $$P(e_{2}^{k+1})(X,b_{2},t)\exp(-\gamma t)\leq E_{k+1},$$ $(\ref{5e4})$ implies  
$$P(e_{3}^{k+2})(X,a_{3},t)\exp(\beta (L_2-S_{2})-\gamma t)\leq\max\{E_{k}\exp(-\beta'_1S_1),E_{k+1}\exp(\beta L_2)\}.$$
Thus
$$P(e_{3}^{k+2})(X,a_{3},t)\exp(-\gamma t)\leq\max\{E_{k}\exp(-\beta (L_{2}-S_{2})-\beta'_1S_1),E_{k+1}\exp(\beta S_2)\}.$$
We choose $\beta=-\beta_2'=-\beta'_1\frac{S_1}{L_{2}}$ such that $-\beta (L_{2}-S_{2})-\beta'_1S_1=\beta S_2$. Then, we have that
\begin{equation}
P(e_{3}^{k+2})(X,a_{3},t)\exp(-\gamma t)\leq\max\{E_{k},E_{k+1}\}\exp(-\beta'_2S_2).
\label{5e8}
\end{equation}
\h Using the same techniques as the ones that we use to achive $(\ref{5e7})$ and $(\ref{5e8})$, we can prove that
\begin{equation}
P(e_{j}^{k+j-1})(X,a_{j},t)\exp(-\gamma t)\leq\max\{E_{k},\dots,E_{k+j-2}\}\exp(-\beta'_{j-1}S_{j-1}), 
\label{5e9}
\end{equation}
where $\beta_j'=\beta_1'\frac{S_1}{L_{2}}\dots\frac{S_{j-1}}{L_{j}}$, $j=\{2,\dots,I-1\}$.
\\ {\h\bf Step 2.3:} Convergence result
\\\h From the fact that $\bar\epsilon=\sqrt \frac{\gamma}{2}\frac{S_1\dots S_{I-1}}{L_2\dots L_{I-1}}$, and
$$\bar E_k=\max_{j\in[0,I-1]}\{E_{k+j}\},$$
$(\ref{5e6})$ and $(\ref{5e9})$ give us
$$\bar E_{k+1}\leq \bar E_k\exp(-\bar\epsilon),\forall k\in\mathbb{N},$$
or 
$$\bar E_{n}\leq \bar E_0\exp(-n\bar\epsilon),\forall n\in\mathbb{N}$$
\\\h Hence $E_k$ tends to $0$ as $k$ tends to infinity. Which gives 
$$\lim_{k\to\infty}\max_{j=\overline{1,I}}||P(e_j^k)||_{C(\overline{\Omega_j\times(0,T)})}=0,$$
and we get the linear convergence of the algorithm.
\section*{Acknowledgements:} The author would like to thank his thesis advisor, Professor Laurence Halpern, for her patience, kindness and her help while training this work.
 
\bibliographystyle{plain}\bibliography{SIAM-DDMforSemilinearHeat}

\begin{thebibliography}{10}

\bibitem{Benamou:1997:DDM}
Jean-David Benamou and Bruno Despr{\`e}s.
\newblock A domain decomposition method for the {H}elmholtz equation and
  related optimal control problems.
\newblock {\em J. Comput. Phys.}, 136(1):68--82, 1997.

\bibitem{Burrage:1996:PPW}
Kevin Burrage, Carolyn Dyke, and Bert Pohl.
\newblock On the performance of parallel waveform relaxations for differential
  systems.
\newblock {\em Appl. Numer. Math.}, 20:39--55, 1996.

\bibitem{Cazenave:1998:AIS}
Thierry Cazenave and Alain Haraux.
\newblock {\em An introduction to semilinear evolution equations}, volume~13 of
  {\em Oxford Lecture Series in Mathematics and its Applications}.
\newblock The Clarendon Press Oxford University Press, New York, 1998.
\newblock Translated from the 1990 French original by Yvan Martel and revised
  by the authors.

\bibitem{Davis:1989:HKS}
E.~B. Davies.
\newblock {\em Heat kernels and spectral theory}, volume~92 of {\em Cambridge
  Tracts in Mathematics}.
\newblock Cambridge University Press, Cambridge, 1989.

\bibitem{Descombes:2010:SWR}
St´ephane Descombes, Victorita Dolean, and Martin~J. Gander.
\newblock {S}chwarz waveform relaxation methods for systems of semi-linear
  reaction-diffusion equations.
\newblock In {\em Domain Decomposition Methods}, 2009.

\bibitem{Evans:1998:PDE}
Lawrence~C. Evans.
\newblock {\em Partial differential equations}, volume~19 of {\em Graduate
  Studies in Mathematics}.
\newblock American Mathematical Society, Providence, RI, 1998.

\bibitem{Friedman:1964:PDE}
Avner Friedman.
\newblock {\em Partial differential equations of parabolic type}.
\newblock Prentice-Hall Inc., Englewood Cliffs, N.J., 1964.

\bibitem{Gander:2007:OSW}
M.~J. Gander and L.~Halpern.
\newblock Optimized {S}chwarz waveform relaxation methods for advection
  reaction diffusion problems.
\newblock {\em SIAM J. Numer. Anal.}, 45(2):666--697 (electronic), 2007.

\bibitem{Gander:1998:OCO}
M.~J. Gander, L.~Halpern, and F.~Nataf.
\newblock Optimal convergence for overlapping and non-overlapping {S}chwarz
  waveform relaxation.
\newblock In {\em Eleventh {I}nternational {C}onference on {D}omain
  {D}ecomposition {M}ethods ({L}ondon, 1998)}, pages 27--36 (electronic).
  DDM.org, Augsburg, 1999.

\bibitem{Gander:1999:WRA}
Martin~J. Gander.
\newblock A waveform relaxation algorithm with overlapping splitting for
  reaction diffusion equations.
\newblock {\em Numer. Linear Algebra Appl.}, 6(2):125--145, 1999.
\newblock Czech-US Workshop in Iterative Methods and Parallel Computing, Part 2
  (Milovy, 1997).

\bibitem{Gander:1999:OSM}
Martin~J. Gander, Laurence Halpern, and Frederic Nataf.
\newblock Optimized {S}chwarz methods.
\newblock In {\em Domain decomposition methods in sciences and engineering
  ({C}hiba, 1999)}, pages 15--27 (electronic). DDM.org, Augsburg, 2001.

\bibitem{Gander:1998:STC}
Martin~J. Gander and Andrew~M. Stuart.
\newblock Space time continuous analysis of waveform relaxation for the heat
  equation.
\newblock {\em SIAM J.}, 19:2014--2031, 1998.

\bibitem{Gander:2002:OSW}
Martin~J. Gander and Hongkai Zhao.
\newblock Overlapping {S}chwarz waveform relaxation for the heat equation in
  {$n$} dimensions.
\newblock {\em BIT}, 42(4):779--795, 2002.

\bibitem{Giladi:2002:STD}
Eldar Giladi and Herbert~B. Keller.
\newblock Space-time domain decomposition for parabolic problems.
\newblock {\em Numer. Math.}, 93(2):279--313, 2002.

\bibitem{Lieberman:1996:SOP}
Gary~M. Lieberman.
\newblock {\em Second order parabolic differential equations}.
\newblock World Scientific Publishing Co. Inc., River Edge, NJ, 1996.

\bibitem{Lions:1987:OSA}
P.-L. Lions.
\newblock On the {S}chwarz alternating method. {I}.
\newblock In {\em First {I}nternational {S}ymposium on {D}omain {D}ecomposition
  {M}ethods for {P}artial {D}ifferential {E}quations ({P}aris, 1987)}, pages
  1--42. SIAM, Philadelphia, PA, 1988.

\bibitem{Lions:1989:OSA}
P.-L. Lions.
\newblock On the {S}chwarz alternating method. {II}. {S}tochastic
  interpretation and order properties.
\newblock In {\em Domain decomposition methods ({L}os {A}ngeles, {CA}, 1988)},
  pages 47--70. SIAM, Philadelphia, PA, 1989.

\bibitem{Lions:1990:OSA}
P.-L. Lions.
\newblock On the {S}chwarz alternating method. {III}.\ {A} variant for
  nonoverlapping subdomains.
\newblock In {\em Third {I}nternational {S}ymposium on {D}omain {D}ecomposition
  {M}ethods for {P}artial {D}ifferential {E}quations ({H}ouston, {TX}, 1989)},
  pages 202--223. SIAM, Philadelphia, PA, 1990.

\bibitem{Lui:2001:OSM}
S.~H. Lui.
\newblock On {S}chwarz methods for monotone elliptic {PDE}s.
\newblock In {\em Domain decomposition methods in sciences and engineering
  ({C}hiba, 1999)}, pages 55--62 (electronic). DDM.org, Augsburg, 2001.

\bibitem{Lui:2003:OMI}
S.-H. Lui.
\newblock On monotone iteration and {S}chwarz methods for nonlinear parabolic
  {PDE}s.
\newblock {\em J. Comput. Appl. Math.}, 161(2):449--468, 2003.

\bibitem{Quittner:2007:SPP}
Pavol Quittner and Philippe Souplet.
\newblock {\em Superlinear parabolic problems}.
\newblock Birkh\"auser Advanced Texts: Basler Lehrb\"ucher. [Birkh\"auser
  Advanced Texts: Basel Textbooks]. Birkh\"auser Verlag, Basel, 2007.
\newblock Blow-up, global existence and steady states.

\end{thebibliography}
\end{document}